\documentclass[12pt,leqno]{amsart}
\usepackage{amsmath}
\usepackage{graphicx,color}
\newtheorem{theorem}{Theorem}[section]

\newtheorem{lemma}[theorem]{Lemma}
\newtheorem{proposition}[theorem]{Proposition}
\theoremstyle{definition}
\newtheorem{definition}[theorem]{Definition}
\theoremstyle{remark}

\numberwithin{equation}{section}

\def\R{\mathbb{R}} 
\def\C{\mathbb{C}} 
\def\D{\mathbb {D}}

\def\Re{{ Re\,}}
\def\Im{{ Im\,}}

\begin{document}

\sloppy

\title{Levi-flat filling of real two-spheres in symplectic manifolds (II)}

\author{Herv{\'e} Gaussier{*} and Alexandre Sukhov{**}}



\date{\number\year-\number\month-\number\day}

\maketitle

{\small
* Universit\'e Joseph Fourier, 100 rue des Maths, 38402 Saint Martin d'H\`eres, France, herve.gaussier@ujf-grenoble.fr

** Universit\'e des Sciences et Technologies de Lille, Laboratoire
Paul Painlev\'e,
U.F.R. de
Math\'e-matique, 59655 Villeneuve d'Ascq, Cedex, France,
 sukhov@math.univ-lille1.fr
}

\bigskip
Abstract. We consider a compact  almost complex manifold $(M,J,\omega)$ with  smooth Levi convex boundary $\partial M$ and a tame symplectic form $\omega$.  Suppose that $S^2$ is  a real two-sphere, containing complex elliptic and hyperbolic points and  generically embedded into $\partial M$. We prove a result on filling $S^2$ by holomorphic discs.

\bigskip
R\'esum\'e. On consid\`ere une vari\'et\'e presque complexe $(M,J,\omega)$ avec la fronti\`ere Levi conv\`exe $\partial M$ et une tame forme sympl\'ectique 
$\omega$. Soit  $S^2$ une 2-sphere r\'eelle  avec des points   elliptiques  et hyperboliques,  plong\'ee g\'en\'eriquement dans  $\partial M$. On d\'emontre un r\'esultat sur le remplissage de  $S^2$ par des disques holomorphes.

\section{Introduction}
This expository paper is the second part  of \cite{ga-su}. We keep the same notations and terminology.  Our main result is the following.
\begin{theorem}
\label{hyperb-theorem}
Let $(M,J,\omega)$ be a compact  almost complex manifold of complex dimension 2 with a
tame symplectic form $\omega$ and smooth boundary $\partial M$. Let also $S^2$
be a real $2$-sphere embedded into $\partial M$. Assume that the following assumptions hold:
\begin{itemize}
\item[(i)] $M$ contains no non-constant  $J$-holomorphic spheres. 
\item[(ii)] the boundary $\partial M$ of $M$ is a smooth  Levi
convex hypersurface containing no non-constant $J$-holomorphic discs.  
\item[(iii)] $S^2$ has only  elliptic and good
hyperbolic complex points. Furthermore, $\partial M$ is strictly Levi convex in a neighborhood of every hyperbolic point.
\end{itemize} 
Then after an arbitrarily small $C^2$-perturbation near hyperbolic points there exists  a smooth  Levi-flat hypersurface $\Gamma \subset M$ with boundary $S^2$. This hypersurface is foliated by $J$-holomorphic discs with boundaries attached to $S^2$.
\end{theorem}

In \cite{ga-su}  we studied the filling of a two-sphere containing only
isolated elliptic complex points. The present work is devoted to the more
general case where hyperbolic points occur. We impose on $M$ and $\partial M$
precisely the same assumptions as in \cite{ga-su}.

Conditions (i)-(ii) are essential to extend a local filling of $S^2$ by
boundaries of pseudoholomorphic discs, starting at an elliptic point, up to
hyperbolic points without appearance of sphere or disc bubbles. Assumption (i) may be weakened (see the discussion in \cite{ga-su}), however (ii) is in general necessary (see the last section). Condition
(iii) is a technical assumption but it is sufficient for consistent applications. A precise definition of a good hyperbolic point is given in the next section. We mention that by definition the almost complex structure $J$ is supposed to be 
integrable near a good hyperbolic  point and each hyperbolic point can be made a good one by an arbitrary small $C^2$ deformation of the sphere $S^2$ near this point. In the last section we include a detailed discussion and comparison of this result with known results obtained by several authors.

We thank the referee for several helpful remarks and suggestions improving the
paper. This work was partially done when the second author was visiting the
Indiana University during the Fall 2011. He thanks this institution for
excellent conditions for work.

\section{Good hyperbolic points: Bedford-Klingenberg's analysis}

Denote by  $S^2_*$  the set of totally real points in $S^2$.
Let $f$ be  a Bishop disc  for $S^2$. We call $f$    {\it hyperbolic} if there exists a finite set of points $\zeta_j \in \partial \D$, $j=1,...,k$, such that every point $f(\zeta_j)$ is a hyperbolic point of $S^2$ and $f(\zeta) \in S^2_*$ for every $\zeta \in \partial \D \backslash \{\zeta_1,\dots,\zeta_k\}$. 
We recall some properties of hyperbolic discs in a neighborhood of a
hyperbolic point. These results were obtained by E.Bedford-W.Klingenberg
\cite{be-kl} in the case where $(M,J)$ coincides with a strictly pseudoconvex
domain in   $(\C^2,J_{st})$. {\it Everywhere in this section we suppose that the almost complex structure
$J$ is integrable near every hyperbolic point of $S^2$.}  This assumption will
be crucially used.

\subsection{Good hyperbolic points.}  Let $S^2$ be a two-sphere generically embedded into an almost complex manifold
$(M,J)$ of complex dimension $2$. In what follows it is convenient to view $M$
as a smoothly bounded subdomain in an almost complex manifold $\tilde M$.

Suppose that  $p$ is a hyperbolic point of
$S^2$.
Since  the almost complex structure $J$ is integrable in a neighborhood of $p$, there exist complex coordinates $(z,w) \in \C^2$ defined in a neighborhood of $p$ such that in these coordinates the origin corresponds to $p$, $J = J_{st}$ in a neighborhood of the origin and $S^2$ is locally  defined by the expression :

\begin{eqnarray} 
\label{hyperb0}
w = z \overline z + \gamma \Re z^2  + o(\vert z \vert^2)
\end{eqnarray}
where $\gamma > 1$.  Then  we can represent  the boundary $\partial M$ near
the origin in the form $\partial M = \{ \rho = 0 \}$ with

\begin{eqnarray*}
\rho(z,w)= -\Re w  + \alpha_1 \Im w  + \vert z \vert ^2 + \alpha_2 \vert w
\vert^2  + \gamma \Re z^2 + \Re [(\alpha_3 + i \alpha_4) z\overline w] +  o(\vert z \vert^2);
\end{eqnarray*}
here $\alpha_j$ are real constants. We will call such a $\rho$ {\it a local
  defining function} of $\partial M$ near a good hyperbolic point.

After an arbitrarily small $C^2$-deformation near the origin  $S^2$ can be transformed to the model quadric :

\begin{eqnarray} 
\label{hyperb1}
w = \psi(z) = z \overline z + \gamma \Re z^2  
\end{eqnarray}
which coincides with the initial sphere outside a neighborhood of $p$. This 
follows by multiplying the $o(|z^2|)$-term in the right-hand side by a suitable cut-off
function vanishing near the origin, see \cite{be-kl}. More precisely, denote
the $o(|z|^2)$-term in (\ref{hyperb0}) by $\phi(z)$. Then we replace $\phi(z)$ by
$\chi(\vert z \vert/\varepsilon)\phi(z)$, where $\chi$ is a smooth function
with $\chi(t) = 0$ for $t < 1$ and $\chi(t) = 1$ for $t > 2$. It is easy to
see that the quantity $(1 - \chi(\vert z \vert/\varepsilon)) \phi(z)$ tends to
$0$ in the $C^2$-norm as $\varepsilon \to 0$.

The boundary $\partial M$ also must be slightly deformed near the origin in order to contain the perturbed  $S^2$ given by (\ref{hyperb1}). Since the hypersurface $\partial M$ is
strictly Levi convex, it remains strictly Levi convex after a small
$C^2$-deformation. Thus, in the present paper we  deal with hyperbolic points
which can be written in the form (\ref{hyperb1}). Throughout the rest of the
paper we keep the notation $\rho$ and $\psi$ for the local
defining functions of $\partial M$ and $S^2$ respectively, near a good hyperbolic point, introduced above.

\bigskip

Now we describe an additional restriction 
 on the values of the parameter $\gamma$ imposed through the present paper; it considerably simplifies the study of hyperbolic discs.
 Consider the proper holomorphic map $H: \C^2 \to \C^2$ given by
$$H(z,w) = (z,zw + \frac{\gamma}{2}(z^2 + w^2)).$$

Then $H$ determines a two-fold branched covering of $\C^2$. The pull-back $H^{-1}(S^2)$ in a neighborhood of the origin consists of two totally real subspaces  $E_1 = \{ w = \overline{z} \}$ and $E_2 = \{ w = -\overline{z} - (2/\gamma)z \}$. Denote by $\tau_j$ the antiholomorphic involution with respect to $E_j$. It is easy to see that 
$$\tau_1 = \tau \circ \left (
\begin{array}{cc}
0&1\\
1& 0
\end{array}
\right )$$
and 
$$\tau_2 = \tau \circ \left (
\begin{array}{cc}
-2/\gamma& -1\\
(2/\gamma)^2 - 1& 2/\gamma
\end{array}
\right )$$
where $\tau$ denotes the usual conjugation in $\C^2$. Since the matrices in the expressions of $\tau_1$ and $\tau_2$ commute with $\tau$, we may with some abuse of terminology view $\tau_j$ as elements of the group $GL(2,\R)$.

\begin{lemma}\label{dense-lemma}
There exists a dense subset $\Lambda$ in $]0,+\infty[$ such that for every $\gamma \in \Lambda$, the involutions $\tau_1$ and $\tau_2$ generate a finite group isomorphic to the dihedral group $D_{2n}$  for some integer $n$ (in general, depending on $\gamma$). 
\end{lemma}
\proof Every matrix 
$$C = \left (
\begin{array}{cc}
c_1&c_2\\
c_2& c_1
\end{array}
\right )$$
commutes with $\tau_1$. Setting $a = 2/\gamma$, we have
\begin{eqnarray*}
& &\tilde\tau_2 = C^{-1} \tau_2 C \\
& &=\frac{1}{c_1^2 - c_2^2}\left (
\begin{array}{cc}
-ac_1^2 - a^2 c_1c_2 - ac_2^2&-c_1^2 - 2ac_1c_2 - (a^2 - 1)c_2^2\\
(a^2-1)c_1^2 + 2c_1c_2a + c_2^2& ac_1^2 + a^2 c_1c_2 + a c_2^2
\end{array}
\right )
\end{eqnarray*}
Since the trace of $\tilde \tau_2$ is zero, it suffices to make it symmetric
in order to assure that this matrix is orthogonal (perhaps, after a multiplication by a
suitable diagonal matrix). The condition for $\tilde\tau_2$ to be symmetric is 
$$ax^2 + 4 x + a = 0$$
with $x = c_1/c_2$. This equation admits a real solution since $0 < a <
2$.The vector $\left (
\begin{array}{c}
\nu\\
1
\end{array}
\right )$ generating the axis of reflection for $\tilde \tau_2$ satisfies 
$$\tilde \tau_2   \left (
\begin{array}{c}
\nu\\
1
\end{array}
\right ) = \left (
\begin{array}{c}
\nu\\
1
\end{array}
\right )$$
From this we deduce
$$\nu = \frac{-c_1^2 - 2ac_1 + 1 - a^2}{(1+a)c_1^2 + a^2 c_1 + a -1}.$$
This implies that $\nu$ is a non-constant algebraic function of the variable
$\gamma = 2/a$. 

The group generated by $\tau_1$ and $\tilde \tau_2$ is generated by orthogonal
reflections about the vectors $\left (
\begin{array}{c}
1\\
1
\end{array}
\right )$ and $\left (
\begin{array}{c}
\nu\\
1
\end{array}
\right )$. If the angle between these vectors is a rational multiple of $\pi$
of the form $n\pi/m$ with relatively prime numbers $n$ and $m$, then $\tilde \tau_2$
and $\tau_1$ generate the dihedral group $D_{2m}$. Since this occurs for a
dense set $1 < \gamma < \infty$, we obtain the conclusion. Q.E.D.

\begin{definition}
\label{def-good}
Suppose that $S^2$ is contained in a strictly Levi convex hypersurface near a hyperbolic point $p \in S^2$. This point  is called {\it good}  if 
\begin{itemize}
\item[(i)] $J$ is integrable in a neighborhood of $p$,
\item[(ii)] there exist local holomorphic coordinates near $p$ such that $S^2$ has the form (\ref{hyperb1}),
\item[(iii)] $\gamma $ satisfies the conclusion of Lemma~\ref{dense-lemma} i.e. $\gamma \in \Lambda$.
\end{itemize}
\end{definition}

\subsection{Analytic extension past a good hyperbolic point.} Now we describe the behaviour of hyperbolic discs near a good hyperbolic point.

\begin{lemma}\label{ext-lemma}
Suppose that  $S^2$ is given by (\ref{hyperb1}) near the origin which is a good hyperbolic point. Let $U$ be a neighborhood of the origin in $\C^2$ and $Y$ be a closed  complex purely  1-dimensional variety in $U \backslash S^2$ (i.e. every irreducible component of $Y$ has complex dimension 1, see \cite{ch1} for more details). Then (possibly  after shrinking $U$) $Y$ is contained in a closed complex purely $1$-dimensional subset $X$ in $U$.
\end{lemma}
Thus, a complex 1-dimensional analytic set $Y$ extends analytically past the
origin. This key result was proved in \cite{be-kl} under an additional
assumption that $Y$ is a graph of a holomorphic function continuous up to the
boundary on a domain in $\C$ whose boundary contains the origin. However the argument still works without this assumption. 
\proof The pull-back 
$V_0 = H^{-1}(Y)$ is a closed complex 1-dimensional  subset  in $H^{-1}(U \backslash
S^2) = H^{-1}(U) \backslash (E_1 \cup E_2)$. Of course, an analytic subset defined outside a single real analytic totally real manifold can be extended through this manifold by the reflection principle \cite{ch1}. Our case is more subtle since the standard reflection principle cannot be applied at the origin where the totally real planes interesect. So we first use the involutions $\sigma_j$ and the reflection principle at points of $E_1 \cup E_2$ off the origin in order to obtain an analytic subset in a punctured neighborhood of the origin.

By Lemma \ref{dense-lemma} the group $G$ generated by the reflections about $E_1$ and $E_2$ is finite. Set $\tilde U = H^{-1}(U)$. Then $\tilde U$ is a neighborhood of the origin and $V := \cup_{\sigma \in G} \sigma(V_0)$ is a complex 1-dimensional subset in 
$\tilde U \backslash (\cup_{\sigma \in G}\sigma(E_1 \cup E_2))$. Since every
$\sigma \in G$ is an antiholomorphic reflection with respect to $E_1$ or
$E_2$,  the reflection principle for complex analytic varieties (see for
instance \cite{ch1}) implies that the closure $\overline{V}$ of $V$ is an
analytic subset in $\tilde{U} \backslash \{ 0 \}$. Then by the Remmert-Stein  removal singularities theorem (see for instance \cite{ch1})  the closure  $\overline{V}$ is a complex analytic variety in $\tilde{U}$. Since $H$ is proper, the image $H(\overline{V})$ is a complex 1-dimensional variety in $U$ containing $Y$.  Q.E.D.

\bigskip

Let $f:\D \to M$ be a non-constant map.Then $f(\D)$ is contained
in $M$, i.e. for every point  $\zeta \in
\D$ its image $f(\zeta)$ is not on  the boundary $\partial M$. This is a
consequence of the $J$-convexity of $\partial M$ and of the assumption
that the boundary $\partial M$  contains no $J$-holomorphic discs, see
\cite{ga-su}. Let $E$ be a non-empty subset of $\partial \D$. The cluster set
$C(f,E)$ of  $f$ on $E$ is defined as the set of all points $p$ such that there
exists a sequence $(\zeta_n)$ in $\D$ converging to a point $\zeta \in E$ with 
$f(\zeta_n) \to p$. A map $f : \D \to M$ is proper if and only if the cluster
set $C(f,\partial \D)$ is contained in $\partial M$ (one can say that in this
sense the boundary values of $f$ belong to $\partial M$). Because of the above
remark  for a proper $J$-holomorphic map  $f:\D \to M$ one has 
 $C(f,\partial \D) = \overline{f(\D)}\backslash f(\D)$. 
Recall the well-known fact that the cluster set $C(f,\partial \D)$
 is connected. Indeed, if not, there are two disjoint open sets $U$ and $V$
 such that $C(f,\partial \D) \subset U \cup V$ and $C(f,\partial \D)$ meets
 each of the sets $U$ and $V$. Then for $r < 1$ close enough to $1$ the connected set
 $f(\{ r < \vert \zeta \vert < 1\})$ is contained in $U \cup V$ and
 intersects each set  $U$ and $V$ which is impossible. 

\bigskip

Let now $f:\D \to M$ be a proper $ J$-holomorphic map such that the cluster set
$C(f,\partial \D)$ is contained in $S^2$. We will see in Section 4 that this
 condition always hold for a disc which is the limit of a sequence of Bishop's
 discs attached to the totally real part $S^2_*$ of the sphere $S^2$. 

Let $p  \in S^2$ be a
 point of $S^2$ and $p \in C(f,\partial \D)$. Consider an open neighborhood
 $U_p$ of $p$. Choosing $U_p$ small enough, we can take it in the form $U_p =
 \{ \phi < 0 \}$ where $\phi$ is a strictly $J$-plurisubharmonic function on $U_p$. It suffices to fix a local coordinate system $(z_1,z_2)$
 centered at $p$ such that in these coordinates $J(0) = J_{st}$ and set
 $\phi(z) = \vert z_1 \vert^2 + \vert z_2 \vert^2 - \varepsilon$ with $\varepsilon > 0$ small enough. Denote by $B(p,\varepsilon) = \{ \phi < 0 \}$ the ball of radius $\varepsilon$ centered at $p$.
This construction may be performed for all points $q$ of $S^2$. It  provides a
family of neighborhoods $(U_q)$, $q \in U_q$, and a family of strictly $J$-plurisubharmonic functions on $U_q$, $(\phi_q)$, depending smoothly on $q$, such that $\phi_p = \phi$.  We use the notation $B(q,\varepsilon) = \{ \phi_q < \varepsilon \}$.

\begin{lemma}
\label{disc-connect}
For $\varepsilon > 0$ small enough each  connected component of  the intersection $f(\D) \cap B(p,\varepsilon)$ is a disc (more precisely, the image of a $J$-holomorphic disc). 
\end{lemma}
\proof   Fix $\varepsilon
> 0$ small enough such that every   ball $B(q,\varepsilon)$ is compactly
contained in the neighborhood $U_q$ of $q$ chosen above, for all $q \in S^2$. Let $G$ be an open
connected component of the pull-back $f^{-1}(B(p,\varepsilon))$. Suppose by
contradiction that $G$ is not simply connected. If the closure $\overline G$ is simply
connected, the subharmonic function $\rho \circ f$ achieves its maximum (equal
to $\varepsilon$) at an interior point of $\overline G$ which gives a
contradiction. Consider the case where $\overline G$ is not simply connected
i.e. $G$ has at least one hole $H$ with non-empty interior. Then $C = f(H)$ is a
compact $J$-holomorphic curve  with boundary $\partial C \subset
B(p,\varepsilon)$. Since $C(f,\partial \D)$ is contained in $S^2$, we can
choose $\varepsilon$ small enough such that the curve $C$ is compactly contained in the
neighborhood $U = \cup_{q \in S^2} U_q$ of $S^2$. On the other hand, the curve
$C$ is not contained in
$B(p,\varepsilon)$. Then, considering the family of balls $B(q,\tau)$ and their
translations along the real normals to $\partial M$, we can find by
continuity suitable $q$ and $\tau > 0$ such that the $J$-holomorphic curve $C$
touches the boundary of the ball $B(q,\tau)$ from inside at some point $a$. But the boundary $\partial B(q,\tau)$ is strictly $J$-convex and  admits
a strictly plurisubharmonic defining function near $a$ as described above. Then the restriction of this function on $C$ is a subharmonic function admitting a local maximum at an interior point. This contradicts the maximum principle. Thus, every open component of  $f^{-1}(B(p,\varepsilon))$ is simply connected. Reparametrizing it by $\D$ via the Riemann mapping theorem, we conclude. Q.E.D.

We point out that the above lemma claims that every single connected component of  the intersection $f(\D) \cap B(p,\varepsilon)$ is a disc. {\it A priori}, such an intersection can have several connected components.

\bigskip

Let $p \in S^2$ be a good hyperbolic point.  Consider   a proper
$J$-holomorphic map  $f:\D \to M$ such that the point $p$  belongs to the
cluster set $C(f,\partial \D) \subset S^2$. Let $U$ be a coordinate neighborhood
of $p$ provided by Definition \ref{def-good}; in particular,  we identify $p$ with the
origin. Then  $Y = f(\D) \cap U$ is a closed complex 1-dimensional subset in
$U \backslash S^2$. By Lemma \ref{disc-connect} this set consists of a finite
number of holomorphic discs in $M \cap U$; we choose one of them and again denote it by $Y$. By Lemma \ref{ext-lemma} the variety $Y$ extends as a complex 1-dimensional set $\tilde Y$  past $S^2$. By the uniqueness theorem for complex analytic sets \cite{ch1} there exists a unique irreducible component of $\tilde Y$ containing $Y$; we still denote this component by $\tilde Y$. In particular, $Y$ and $\tilde Y$ do not contain  the $w$-axis (if not, $\tilde Y$ would coincide with the $w$-axis by the uniqueness theorem).

 Let now  $\pi:(z,w) \mapsto z$ be the canonical projection. Since $\tilde Y$ does not contain the $w$-axis, the restriction $\pi\vert_{\tilde Y}$ is proper when $U$ is small enough. More precisely, the intersection  $\pi^{-1}(0)  \cap \tilde Y$ is discrete near the origin and taking a neighborhood $U = U' \times U''$ small enough, $U'$ and $U''$ being neighborhoods of the origin in $\C$, we obtain that $\pi: \tilde Y \cap U \to U'$ is proper (see \cite{ch1}).

Then $Y$ is the graph $\{ w = g(z) \}$ of a function $g$ holomorphic in a
domain $D$ in $\C$, and the boundary of $D$ contains the origin. Since the
variety $Y$ extends analytically past $S^2$, the function $g$ is continuous up
to the boundary of $D$ and extends past the origin as a multivalued complex
analytic function. Furthermore,  $g$ has a Puiseux expansion at the
origin. Since $g(z) = \psi(z)$ on $\partial D$, we have $g(z) = O(z^2)$ for $z
\in \partial D$ and  all terms of this expansion are $O(z^2)$.
Then there exists a positive integer $m$ such that in a neighborhood of the origin $g$ is given in $D$ by the Puiseux expansion
\begin{eqnarray}
\label{hyperb2}
g(z) = \sum_{k \geq 2m} g_k z^{k/m}.
\end{eqnarray}

The representation (\ref{hyperb2}) gives much useful information about the
behaviour of a hyperbolic disc near a good hyperbolic point. For instance, we obtain that  the boundary of $f(\D)$ is a
continuous curve (with a finite number of real analytic components
intersecting at the origin) near the origin. Now under the assumption that the area of
$f$ is bounded (which always holds in our situation) one can easily show that the map
$f$ itself is necessarily continuous on the closed disc $\overline\D$. The proof is based on a classical argument form the one-variable theory of conformal maps. This justifies the terminology ``a hyperbolic disc'' since by
definition such a disc  is a Bishop disc and so is continuous up to the
boundary. We postpone the proof of this fact  to Section 4, but in the rest of this
section one may assume that a hyperbolic disc is continuous on $\overline\D$.

\vskip 0,1cm
\subsection{Approaching a good hyperbolic point by a disc.} Following \cite{be-kl} we study a local behaviour of  hyperbolic discs near a good hyperbolic point. Since these results will be crucially used, for reader's convenience we include the proofs.

We begin with  a more precise information about the Puiseux expansion (\ref{hyperb2}):

\begin{lemma}
\label{puiseux-lemma}
In the expansion (\ref{hyperb2}) the sum is taken over the set of integers $k$ satisfying $k > 2m$.
\end{lemma}
\proof Suppose by contradiction that $g_{2m}$ in (\ref{hyperb2}) does not vanish. Then 

\begin{eqnarray}
\label{hyp111}
g(z) = a z^{2} + O(z^{2 + 1/m}), 
\end{eqnarray}
with $a \neq 0$. Recall that $g$ is holomorphic on a domain $D$ in $\C$ whose
boundary contains the origin. Denote by $\D^+$ the upper half-disc $\D^+ = \{
\zeta \in \C: \vert \zeta \vert < 1, \Im \zeta > 0 \}$. It is convenient to
assume that the hyperbolic disc $f$ under consideration is defined and
holomorphic on $\D^+$, continuous on $\overline\D^+$ and $f([-1,1]) \subset S^2$
i.e. $f(\D^+)$ is glued to the sphere $S^2$ along the real segment $[-1,1] =
\{ \zeta: -1 \leq \Re \zeta \leq 1, \Im \zeta = 0 \}$; furthermore, $f(0) =
0$. One can always achieve these conditions, reparametrizing $f$ by a suitable
conformal isomorphism. The image $f(\D^+)$ coincides with the graph of $g$ over the domain $D$.  As above, denote by $\pi: \C^2 \to \C$, $\pi(z,w) = z$, the canonical projection. If $\varepsilon > 0$ is small enough, the intersection $\partial D \cap \varepsilon \D$ consists of two real curves $\gamma_+ = (\pi \circ f)([-1,0])$ and $\gamma_- = (\pi \circ f)([0,1])$. Since $S^2$ is given by (\ref{hyperb1}) near the origin, 
$\gamma_+ \cup \gamma_-$ is contained in $\{ z: \Im g(z) = 0 \}$. Hence
$\gamma_\pm$ are real analytic 1-dimensional sets. By (\ref{hyp111}) each
$\gamma_\pm$ is tangential at the origin to one of the real lines satisfying
the equation $\{ \Im a z^2 = 0 \}$. Therefore, the domain $D$ at the origin is
asymptotic to an angle of size $\kappa$ (i.e. the curves $\gamma_\pm$ are
tangent at the origin to the rays bounding this angle) which is a non-zero
integer multiple of $\pi/2$. Consider separately the  possible cases.

{\it Case 1: $\kappa \geq \pi$}. Recall that in a neighborhood $U$ of the
origin the disc $f(\D^+)$ is contained in the domain $M \cap U = \{ \rho < 0
\}$ where a strictly plurisubharmonic function $\rho$ is a local defining
function for $M$ introduced in Subsection 2.1. The composition  $\phi(z) = \rho(z,g(z))$ is a negative subharmonic function on $D$ and its gradient vanishes at the origin. By assumption on $\kappa$,one can find an open  disc $G$ contained in $D$ such that $\partial G \cap \partial D = \{ 0 \}$. Applying to $\phi$ the Hopf lemma on $G$, we obtain that $\partial\phi(0)/\partial z \neq 0$ which is a contradiction.

{\it Case 2: $\kappa = \pi/2$}. Since $\Im g(z) = 0$ for $z \in \gamma_+ \cup \gamma_-$, we conclude that
$$g(z) = \pm \vert a z^2 \vert + O(\vert z \vert^{2 + 1/m}), z \in \gamma_+ \cup \gamma_-$$
and furthermore, $g$ has opposite signs on $\gamma_+$ and $\gamma_-$. Every curve $\gamma_\pm$ is tangent at the origin to a line forming an angle $\mu_\pm$ with the axis $x$, so that $\mu_- = \mu_+ + \pi/2$. The function $\psi$ (defined in (\ref{hyperb1}) in the polar coordinates $z = re^{i\theta})$ has the form 
 $\psi(z) = r^2\chi(\theta)$ with $\chi(\theta) = 1 + \gamma \cos(2\theta)$. 
Consider a sequence $z_k = r_ke^{i\theta_k}$ in $\gamma_+$ converging to
$0$. Then for every $k$ one has $\psi(e^{i\theta_k}) = r_k^{-2}
g(r_ke^{i\theta_k})$. Passing to the limit as $k \to \infty$, we obtain that
$\psi(e^{i\mu_+})$ is equal to  $\alpha$ or $-\alpha$. Repeating the same argument for $\gamma_-$, we obtain 
$\psi(e^{i\mu_+}) = -\psi(e^{\mu_+ + \pi/2})$. But this means that $\mu_+$ satisfies the equation 
$$1 + \gamma\cos(2\mu_+) = -1 - \gamma\cos(2(\mu_+ + \pi/2)).$$
But this equation does not admit any solution.
Thus, none of the above cases can occur. This implies $a = 0$. Q.E.D.

\vskip 0,2cm
Let $D^{+}$ and $D^-$ denote the connected components of  $\{ z \in \C:
z\overline{z} + \gamma \Re z^2 < 0 \}$. They are  domains defined by the inequalities of type $\{ \vert \Re z \vert < C \vert \Im z \vert \}$.
The sign of $D^{\pm}$ is chosen so
that the domain $D^+$ intersects the axis $\{ Im \ z > 0 \}$. 

\begin{definition}
\label{GoodApproach}
A hyperbolic disc $f$ has a {\it  good (or admissible) approach} at $(0,0)$ if 
\begin{itemize}
\item[(i)] There exists a domain $D \subset \C$ which is asymptotic at $0$ to one of the domains $D^{\pm}$ (that is the boundary $\partial D$ near $0 \in \partial D$ 
is formed by two curves tangent at the origin to the boundary lines of $D^{\pm}$).
\item[(ii)] There exists a  function $g$ holomorphic on $D$ such that in a neighborhood of the origin the image $f(\D)$ is the graph of $g$ over $D$.
\item[(iii)] $g \in C^2(\overline D)$ and $g(0) = g'(0) = g''(0) = 0$.
\end{itemize}
\end{definition}

The first important  consequence of Lemma \ref{puiseux-lemma} is the following assertion:
 \begin{lemma}
 \label{GoodApproach-lemma}
 A hyperbolic disc has a good approach at a good hyperbolic point.
 \end{lemma}
\proof  As above, the function $\psi$ in the left-hand side of (\ref{hyperb1}) has the form $\psi(z) = r^2\chi(\theta)$ with $\chi(\theta) = 1 + \gamma \cos(2\theta)$. One readily sees that there is an angle $0 < \mu < \pi/4$ such that $\{ \chi = 0 \} = \{ \theta = \pi/2 \pm \mu \}$. Let us write the two first terms of the expansion (\ref{hyperb2}):
\begin{eqnarray}
\label{hyp004}
g(z) = a z^{k/m} + b z^{k'/m} + o(z^{k'/m}), 2m < k < k'.
\end{eqnarray}
Consider a point $z = re^{i\theta(r)}  \in \partial D$. It follows from (\ref{hyp004}) and from the identity $g(z) = r^2 \chi(\theta(r))$, $z \in \partial D$, that $\chi(\theta(r)) = O(r^{k/m-2})$. Since the right-hand side converges to $0$ as $r \to 0$, we conclude that $\theta(r)$ converges to  $\pi/2- \mu$ or to $\pi/2 + \mu$ as $r \to 0$ when $z$ belongs to $\gamma_+$ or $\gamma_-$ respectively.  Thus, $\partial D$ is tangent to the lines $\{ \theta = \pi/2 \pm \mu \}$ at the origin. Q.E.D.

Another important consequence is the uniqueness principle for hyperbolic discs. 

\begin{lemma}
 \label{unique-lemma}
 Suppose that the origin is a good hyperbolic point. Let $f_1$ and $f_2$ be
 hyperbolic discs in a neighborhood $U$ of  the origin which are the graphs of
 functions $g_j$ over domains $D_j$. Suppose that one of the following
 conditions holds :
\begin{itemize}
\item[(i)] $D_1 = D_2$,
\item[(ii)] $D_1 \subset D_2$ and $\Re g_1(z) \leq \Re g_2(z)$, $z \in U
  \cap D_1$. 
\end{itemize}
Then $f_1 = f_2$.
 \end{lemma}
 
 \proof (i) This is immediate. Indeed,  since $f_j(\partial \D) \subset S^2$, we have $g_j(z) =
\psi(z)$ for $z \in \partial D$
near the origin i.e. $g_1(z) = g_2(z)$ on an arc of $\partial D$ of positive
length. Then $g_1 \equiv g_2$ by the classical boundary uniqueness theorem for
holomorphic functions.

(ii) First we study more precisely the asymptotic behaviour of $D$ at the
 origin. We use the equality  $\{ \theta = \pi/2 \pm \mu \} = \{ \Im a z^{k/m}
 = 0 \}$ established above. Let $a = r_0 e^{i\theta_0}$. Then the set $\{ \Im
 a z^{k/m} = 0 \}$ consists of the rays $R_l:= \{\theta = (m/k)l\pi - (m/k)\theta_0\}$, $l$ being an integer. For some $l$ the rays $R_l$  and $R_{l+1}$ coincide with $\theta = \pi/2 \pm \mu$ which implies that 
 $\pi m/k = 2 \mu$ and $\theta_0 = \pi l + 2 \pi - 2(k/m)\pi$.  Set 
 $\mu = -\chi'(\pi/2 - \mu) = \chi(\pi/2 + \mu)$. We have $\theta(r) = \pi/2 \pm \mu + o(r^{k/m - 2})$ when $z \to 0$ along $\gamma_\pm$. Expanding the identity 
 $\theta(r) = \chi^{-1}(r^{-2}g(z))$ at the origin we obtain 
 \begin{eqnarray}
 \label{hyp007}
 \theta(r) = \pi/2 \pm \mu - \frac{\vert a \vert}{\nu}r^{k/m-2} + O(r^{(k+1)/m - 2}).
 \end{eqnarray}
 
 It is worth noticing here that the constants $\mu$, $k/m$, $\nu$ are determined by $\gamma$ and are independent of $g$.
 
 Now we proceed the second step of the proof. Consider the Puiseux expansion $g_j(z) = a_j z^{k/m} +...$ for 
 $j = 1,2$. The boundary of every domain $D_j$ is given by (\ref{hyp007}) near the origin. Since $D_1$ is contained in $D_2$ we obtain that $\vert a_1 \vert = \vert a_2 \vert$. For the same reason the rays $R_l$ and $R_{l+1}$, corresponding to $j=1,2$, are defined by the same $l$. Hence $a_1 = a_2$. Next, consider the difference
 $$g_2(z) - g_1(z) = \sum_{n > k} c_nz^{n/m}$$
 Let $c_q$ be the first non-vanishing coefficient. Then 
 $$0 < \Re g_2(z) - \Re g_1(z) =  \Re (c_q z^{q/m}) + o(z^{q/m}).$$
 on $D_1$ near the origin. But $\Re (c_q z^{q/m})$ can be positive only on a sector of angle $\pi m/q$ 
 while $\partial D_1$ is asymptotic to a sector of larger angle $\pi k/m$ : a contradiction. Q.E.D.

\section{Indices }

In the first part of our work \cite{ga-su} we saw that an extension of a 1-parameter family of Bishop discs is determined by some topological characteristic of these discs called in \cite{ga-su} the winding number. In the present work we need to study the behaviour of this invariant when a family of discs extends past a hyperbolic point. For our applications it is  convenient to discuss this notion more conceptually and in full generality. This is the goal of the present section. 

We  recall some known facts concerning the Maslov index and the topological properties of real surfaces in (almost) complex manifolds. 

\subsection{Totally real case : the Maslov index} There are several
possibilities to introduce V.Arnold's conception of  the Maslov index \cite{ar}; we follow \cite{mc-sa, ho-li-si}. Denote by $ S^1$ the unit circle. Let 
$$R(n) = GL(n,\C)/GL(n,\R)$$
be the manifold of totally real n-dimensional subspaces of $\C^n$. Consider the map 
$\kappa: R(n) \to S^1$ defined by
$$\kappa(B \cdot GL (n,\R) ) = \frac{\det (B^2)}{\det ({B}^* B)}, B \in GL(n,\C)$$
where the index star denotes the matrix transposition and the complex conjugation. Let $\gamma: S^1 \to R(n)$ be a continuous map i.e. a loop in $R(n)$. {\it The Maslov index} of $\gamma$ is defined by 
$$\mu(\gamma) = \deg(\kappa \circ \gamma)$$
where $\deg$ denotes the topological degree of a map. V.Arnold \cite{ar}
proved that  two loops in $R(n)$ are homotopic if and only if they have the same Maslov index. According to classical results a complex vector bundle $L$ over  the unit circle $S^1$ is trivial because $\pi_0(GL(n,\C)) = 0$.  If $F$ is a totally real subbundle (with fibers of real dimension $n$) of the trivial bundle $L = S^1 \times \C^n$ , we consider the loop $\gamma: S^1 \to R(n)$ defined by 
\begin{eqnarray}
\label{loop}
\gamma(\zeta) = F_\zeta, \zeta \in S^1.
\end{eqnarray}
Here $F_\zeta$ denotes the fiber of $F$ at the point $\gamma(\zeta)$. Hence it is a totally real subspace of $\C^n$ and can be viewed as an element of $R(n)$. Since $\pi_0(GL(n,\R)) = {\bf Z}_2$, there are two rank $n$ real vector bundles over the circle. Two totally real subbundles are isomorphic as real bundles if their Maslov classes have the same parity.

Let $L \to \D$ be a complex rank $n$ vector bundle over the unit disc $\D$ and $F \subset L \vert_{\partial\D}$ be a totally real subbundle over $\partial\D$. {\it The Maslov index } $\mu(L,F)$ of the pair $(L,F)$ is an integer which can be defined axiomatically by the  properties of isomorphism, direct sum, normalisation and decomposition. Let us describe them.

{\it Isomorphism. } If $\Phi: L_1 \to L_2$ is a vector bundle isomorphism then 
$$\mu(L_1,F_1) = \mu(L_2,\Phi(F_1)).$$

{\it Direct sum.} 
$$\mu(L_1 \oplus L_2, F_1 \oplus F_2) = \mu(L_1,F_1) + \mu(L_2,F_2).$$

{\it Normalisation.} If $L = \D \times \C$ is the trivial line bundle and $F_\zeta = e^{ik\theta/2}\R$, $\zeta = e^{i\theta} \in \partial\D$, then
$$\mu(L,F) = k.$$

We do not describe the decomposition property (see \cite{mc-sa}) since we do not need it here.

In general one can show that if $L = \D \times \C^n$ is the trivial bundle, then 
$$\mu(L,F) = \mu(\gamma)$$
where the  loop $\gamma$ in $R(n)$ is defined by (\ref{loop}). One can use this property as the definition of the Maslov class and then 
 show (see \cite{mc-sa}) that it is independent of a choice of trivialisation. The following trivialisation presented in \cite{mc-sa} is particularly useful.

\begin{lemma} 
\label{trivial}
\label{LemMaslov}For every complex line bundle $L$ over $\D$ and every totally real subbundle $F \subset L\vert_{\partial\D}$ there exists a trivialisation such that in the corresponding coordinates one has 
$F_\zeta = e^{ik\theta/2}$, $\zeta = e^{i\theta} \in \partial \D$.
\end{lemma}

Let now $f$ be a Bishop disc with boundary attached to a totally real manifold
$E$ in an almost complex manifold $M$. The pull-back $f^*(TM)$ is a complex
vector bundle over $\D$ and $f^*(TE)$ is its totally real subbundle over
$\partial \D$. Then the Maslov class $\mu(f^*TM,f^*TE)$ is defined and is
called {\it the Maslov index of the disc} $f$. We denote it   by $\mu_E(f)$.

It was shown in \cite{ga-su} that every elliptic point of the
sphere $S^2$ generates a 1-parameter family of Bishop discs $(f^t)$ whose
boundaries foliate a punctured neighborhood of this point in $S^2$. We
introduced in \cite{ga-su} a topological invariant of $f^t$ called the winding
number. It is equal to $0$ for every disc $f^t$.
Comparing that definition of the winding number 
for a Bishop disc $f^t$ from \cite{ga-su} with Lemma \ref{LemMaslov} one
readily sees that  the Maslov index of $f^t$ also is equal to $0$. Thus, the Maslov index of every disc $f^t$ generated by an elliptic point is equal to $0$. We sum up this in the following assertion.

\begin{lemma}
Let $(f^t)$ be a 1-parameter family of Bishop discs generated by an elliptic point as in \cite{ga-su}. Then the winding number of each disc coincides with its Maslov index. In particular
$$\mu_{S^2}(f^t) = 0.$$
\end{lemma}

\subsection{Index of a complex point} Here we follow \cite{fo,ne}. Let $S$ be
a real surface embedded into an almost complex manifold $M$ of complex dimension
2. Recall that  we assume everywhere that the complex points of $S$ are
isolated and are either elliptic or hyperbolic. Assume for simplicity that the
almost complex structure $J$ is integrable near the complex points of $S^2$;
this assumption can be easily dropped. Let $p$ be a complex point in $S$. There
exist local complex coordinates $(z,w)$ centered at $p$ such that $S$ is
locally the graph  $\{ (z,w): w=g(z) \}$ of   a smooth complex valued function
$g$ defined in a neighborhood of the origin in $\C$ and  such that $g(z) =
O(\vert z \vert^2)$.
The points where $\partial g/\partial \overline{z} \neq 0$ are totally
real. Since the origin is an isolated  complex point,  $\partial g/\partial \overline{z}$ does not vanish
elsewhere in a neighborhood of the origin. The index of $p$, denoted by
$I(p,S)$ is defined as the winding number of the function $\partial g
/ \partial \overline{z}$ around the origin. In the general almost complex situation consider first the case of an elliptic point. Choose local coordinates as in \cite{ga-su}; in particular, 
 $J(0) = J_{st}$. Then  the origin is also an isolated complex point for the standard structure and we can   apply the
above construction using the standard operator  $\partial g
/ \partial \overline{z}$. Since we consider here the case where an almost complex structure is integrable near each hyperbolic point, the above construction can be applied directly there.
We point out that $I(p,S)$ does not depend on a choice of $g$ defining $S^2$ locally. For a totally real point $p \in S$ we set $I(p,S)=0$. One readily sees from this definition that the index of an elliptic point is equal to $+1$ and the index of a hyperbolic point is equal to $-1$. This definition does not rely on an orientation of $S$. 

If $S$ is orientable (which is the case considered in the present paper) the orientation must be taken into account.

\begin{definition} 
A complex point $p$ of $S^2$ is called {\it positive} (resp. {\it negative}) if the orientation of  the tangent space $T_pS^2$ coincides with (resp. is opposite to) the orientation induced by $(M,J)$. 
\end{definition}
We define  $I_+( S)$ (resp. $I_-( S)$) as  the sum of the indices over positive (resp. negative)  complex points in $S$. It is known that $I_+(S) = I_-(S)$ for every closed oriented immersed surface in $\C^2$, see \cite{fo}.

The index of $S$ is defined by $I(S):=\sum_{p \in  S}I(p, S)$; in particular $I( S) = I_+( S)+I_-( S)$. We note that for a sphere $S$ generically embedded into a (almost) complex manifold we have 
$$I(S) = \chi(S) = 2$$
where $\chi(S)$ is the Euler number.

We denote by 
  $e_+( S)$ (resp. $e_-( S)$) the number of positive (resp. negative) elliptic points of $ S$ and by $h_+( S)$ (resp. $h_-( S)$) the number of positive (resp. negative) hyperbolic points. Then from this  definition we get $I_{\pm}( S)=e_\pm( S) - h_\pm( S)$. Moreover, if $S$ is an oriented  real surface embedded to a complex surface $M$ we have :

\begin{lemma}\label{lai-lemma}
\begin{equation}\label{lai-equation}
I_\pm(S) = (1/2)(I(S) \pm c(S)).
\end{equation}
\end{lemma}

Here $c(S)$ denotes the value of the first Chern class $c_1(M)$ on $S$. For example, $c(S)$  always vanishes  when $M = \C^2$; in what follows we assume that this condition always holds  in the present   paper.
In particular we have : $I_\pm(S^2) = 1$ for a two-sphere $S^2$ embedded into an almost complex manifold.

In \cite{fo} a slightly different version of the notion of the index of a loop is used. Let us recall it since this is useful in index computations. The simplest way is to define it in a special system of coordinates. Let $U$ be an open set in the complex plane, $g:U \to \C$ be a smooth complex valued function on $U$, and $E \subset \C^2$ be the graph of $g$:

\begin{eqnarray}
\label{graph}
E = \{ (z,g(z)) \in \C^2: z \in U \}.
\end{eqnarray}

Let now $\gamma = (\gamma_1,\gamma_2): S^1 \to E$ be a loop contained in the set of totally real points of $E$. Then the index $I_E(\gamma)$ is equal to the winding number of the function 
$\theta \in S^1 \mapsto \partial g / \partial \overline z (\gamma_1(\theta)) \in \C
\backslash \{ 0 \}$. In some situations it is not convenient to use a
coordinate  representation (\ref{graph}). If $E$ is embedded to $\C^2$, the
simplest way to compute the index $I_E(\gamma)$ is the following (see
\cite{fo}). Choose continuous vector fields $X_j:S^1 \to \C^2$, $j = 1,2$ such
that for every $\theta \in S^1$ the vectors $X_j(\theta)$ form a basis of
$T_{\gamma(\theta)}E$. Such vector fields exist when $E$ is orientable along
$\gamma$ which is always our case; the method can be extended to the
non-orientable case, see  \cite{fo}. Then the index $I_E(\gamma)$ is equal to
the winding number  (around the origin) of the determinant $\det(X_1(\theta),X_2(\theta))$. This determinant does not vanish on the unit circle since $E$ is totally real. Hence  $I_E(\gamma)$ is correctly defined and is independent on a choice of $X_j$ , in particular, is independent on an orientation of $E$ (but depends, of course, on an orientation of $\gamma$). Indeed, if $X_j'$, $j=1,2$ is another couple of vector fields with similar properties, then $X_j' = AX_j$ for some continuous map $A: S^1 \to GL(2,\R)$. Hence the winding number of $\det(X_1',X_2')$ is equal to the sum of winding numbers of $\det A$ and $\det(X_1,X_2)$. But $\det A$ is a real-valued non-vanishing function on $S^1$, so its winding number is equal to $0$. 
One can give an intrinsic definition of $I_E(\gamma)$
for an oriented  real surface embedded or immersed into a (almost) complex
manifold. This general  definition is  similar to the notion of the
Maslov index so we drop it, see details in \cite{fo}.

We recall here that the Maslov index of the boundary $f\vert_{\partial\D}$ of
a Bishop disc generated by an elliptic point is equal to $0$. On the other hand, a direct computation shows that $I_{S^2}(f\vert_{\partial\D}) = 1$.
It is easy to
see from the above definitions that for a   Bishop disc $f$ the equality  $\mu_{S^2}(f) = 0$ holds if and only if $I_{S^2}(f\vert_{\partial\D}) = 1$. Indeed, choose  canonical coordinates along $f(\D)$ as in \cite{ga-su}; they provide a  trivialisation of bundles from Lemma \ref{trivial}.  Comparing the Maslov index and  the winding number of the determinant as described above (with an obvious choice of vector fields $X_j$), we conclude.

The following useful statement is contained in \cite{fo}.

\begin{proposition}
\label{ProInd}
Let $E$ be an oriented real surface with totally real boundary $\partial E$ . Then 
$$I_+(E) - I_-(E) = I_E(\partial E).$$
\end{proposition}

The above mentioned results are often proved for real surfaces in complex
manifolds. It is easy to see that the proofs remain true
without any changes in the almost complex case since they use only standard differential
geometry and the topological properties of complex vector bundles.

 \section{Approaching good hyperbolic points by families of discs} 

 Let $f^t$,
 $t \in ]0,1[$, be a one-parameter family of embedded $J$-holomorphic discs of Maslov
 index $0$ attached  to the totally real part $S^2_*$ of the sphere
 $S^2$. They have a uniformly bounded area  and  their boundaries
 $f^t(\partial\D)$ foliate an open subset $E$ of $S^2_*$, see \cite{ga-su}. Let
 $(f_k)_k$ be a sequence of such discs corresponding to the values $(t_k)$ of
 the parameter i.e. $f_k = f^{t_k}$. Since the areas are bounded, Gromov's compactness theorem can
 be applied. The case where these discs are separated from the set of complex
 points of $S^2$ is considered in \cite{ga-su}. In that case there are no
 bubbles and after a suitable reparametrization by Mobius transformations the
 sequence $(f_k)$ converges in every $C^m(\overline\D)$ norm to a non-constant Bishop
 disc attached to $S^2_*$. Let us consider now the case where the limit of a sequence of
 discs  touches a good hyperbolic point $p$. The usual version of 
Gromov's compactness theorem deals with totally real manifolds. Here we adapt it to our situation.
Previously we often identified a map $f_k:\D \to M$  and its  image $f_k(\D)$ using 
the same terminology "disc" for both of them. In this section we proceed more carefully and distinguish  the convergence of $(f_k(\D))_k$ as a sequence of sets 
(i.e. non-parametrized holomorphic curves) and the convergence of  $(f_k)_k$ as a sequence of maps (i.e. parametrized holomorphic curves).

\subsection{ Convergence of $(f_k)$ as sets.}  Let $U$ be a neighborhood of a good hyperbolic point $p$. Then $Y_k = f_k(\D) \cap U$ is a closed complex purely 1-dimensional 
analytic subset in $U \backslash S^2$ for every integer $k \geq 1$. Since the areas
of the sets $Y_k$ are uniformly bounded, Bishop's convergence theorem
\cite{ch1} implies (after extracting a subsequence) that the sequence $(Y_k)$ converges on every compact subset of $U \backslash S^2$ to a complex purely 1-dimensional analytic subset $Y$ in $U \backslash S^2$. Here the convergence is in the sense of the Hausdorff distance. Applying Lemma \ref{ext-lemma} to $Y$  we conclude that $Y$ extends to a neighborhood of $p$ as a complex 1-dimensional analytic set. Thus, the convergence of $(f_k(\D))$ as sets near a hyperbolic point is quite simple. Unfortunately, this argument is not precise enough. In order to obtain a more detailed information about the limit analytic set $Y$, we need to study the convergence of the sequence of maps $(f_k)$.

\subsection{Interior convergence of $(f_k)$ as maps}  We proceed in several steps. First, recall some basic notions related to the Gromov compactness theorem that we apply to our sequence $(f_k)_k$ (see \cite{ga-su} and references there for precise definitions and statements). After a suitable reparametrization we may assume that $(f_k)$ converges to a non-constant map $f_\infty$ attached to $S^2$. The convergence is in {\it Gromov's sense} which we briefly discuss.

We say that {\it a spherical bubble} arises at an interior point $q$ of the unit disc $\D$ if one can find a sequence of biholomorphic maps $(\phi_k)$ which "blow up" a neighborhood of $q$ such that the sequence $(f_k \circ \phi_k)$ converges to a non-constant $J$-holomorphic map $g:\C \to M$. Since the areas are bounded, $g$ extends to the whole Riemann sphere as a $J$-holomorphic map and is called a spherical bubble. Similarly we may define {\it disc-bubbles} that can occur only at boundary points of $\D$. Since the areas of $f_k$ are uniformly bounded with respect to $k$, bubbles can occur only in a finite number of points $\Sigma \subset \overline\D$. The sequence $(f_k)$ converges to the limit disc $f_\infty$ uniformly on every compact subset
of $\overline\D \backslash \Sigma$ and  at every point of the set $\Sigma$ a bubble
necessarily arises. After a suitable reparametrization by Mobius
transformations the sequence converges to a prescribed bubble. The images
$f_k(\overline\D)$ converge in the Hausdorff distance to the union of
$f_\infty(\D)$ and of the images of the bubbles. Furthermore, this union is { \it a
  nodal curve}, in particular, is connected. Recall that {\it the standard node}
is the complex analytic set $\{ (z_1,z_2) \in \C^2: z_1 z_2 = 0 \}$. A point
on a complex curve is called a {\it nodal point} if it has a  neighborhood
biholomorphic to the standard node. A boundary nodal point may be defined
using the double of a given complex curve with boundary (namely its
compactification as a Riemann surface, see, for instance, \cite{pa}). A nodal
curve is a compact complex curve with boundary and with a finite number of
interior and boundary nodes (see \cite{iv-sh,pa} for more details). Notice
that in the case under consideration the following simplification
occurs. According to Section 2, every hyperbolic disc extends as a complex
analytic set through a good hyperbolic point. Hence, if a disc-bubble arises
at a good hyperbolic point $p$, this point $p$ will be an interior nodal point for a complex analytic set in a neighborhood of $p$ which  extends "the principal disc" and a bubble. 

In what follows we need a more detailed description of the behaviour of
the sequence $(f_k)$ near a point where a disc bubble arises. We follow
\cite{pa,si,ye}. Let $q$ be a point in $\Sigma \in \partial \D$. Then one can
find a sequence of discs $D_k:= q_k + r_k\D$ of radius $r_k \to 0$, centered at
$q_k \in \partial \D \to q$, and conformal injective maps $\psi_k: D_k \cap \D
\to H:= \{ \zeta \in \C: \Im \zeta > 0 \}$ such that $\psi_k(D_k
\cap \partial \D) \subset \partial H$, $\psi_k(D_k \cap \D) \to H$ and $f_k
\circ \psi^{-1}_k$ converges smoothly on compacts to a disc-bubble $g$
(more precisely, $g$ becomes a disc bubble after a composition with a
conformal isomorphism from $\D$ to $H$). Moreover, $f_k(\partial D_k \cap \D)$
converges to a cusp (nodal) point $f_\infty(q)$ where $f_\infty$ and $g$ intersect each other. 

The connectedness of the limit nodal curve  is often useful to study its global topological properties (the Maslov indices of discs contained in such a curve, their homology classes, etc.). Furthermore, for large $k$ the homology class  $[f_k\vert_{\partial\D}]$ is equal to the sum of $[f_\infty\vert_{\partial\D}]$ and of the homological classes of bubbles boundaries on $S^2$ (see, for instance, \cite{pa,ye}).

\bigskip

Recall that the first step in the proof of Gromov's compactness theorem
\cite{mc-sa,si} considers holomorphic discs with free boundary (i.e. without totally real
boundary conditions) and claims that  after a suitable reparametrization the sequence $(f_k)$ converges to a non-constant disc $f_\infty$ outside a finite set $\Sigma$ as explained above.
By assumption,   $M$ contains no  non-constant $J$-holomorphic
spheres. Since spherical bubbles are non-constant, there are no spherical bubbles. Hence the only possibility is that 
the set $\Sigma$ where bubbles arise is contained in $\partial\D$ and the bubbles must be disc-bubbles. In
particular, the sequence  $(f_k)$ converges to $f_\infty$ uniformly on compact subsets of
$\D$. But then, since the areas of the discs are uniformly bounded, the sequence
$(f_k)$ satisfies the assumptions of the theorem in  \cite{la} which implies that
the cluster set $C(f_\infty,\partial \D) = \overline{f_\infty(\D)}
\backslash f_\infty(\D)$ is contained in $S^2$. The same  holds for the 
boundary disc-bubbles since after suitable reparametrization the sequence $(f_k)$ converges
uniformly on compact subsets of $\D$ to such a bubble.

Our next goal is to study the behaviour of boundary disc bubbles. We will see
that they do not arise at totally real boundary points of a hyperbolic disc
and will give a precise description of bubbling near a good hyperbolic
point. We begin with a useful technical statement mentioned above in Section 2.

\subsection{Boundary continuity of discs forming the  limit nodal curve} Let
$f:\D \to M$ be a non-constant $J$-holomorphic disc from the limit nodal curve (i.e. after
a suitable  reparametrization by a sequence of Mobius transformations the
sequence $(f_k)$ converges to $f$). As we noticed above, the cluster set
$C(f,\partial \D) = \overline{f(\D)}\backslash f(\D)$ of $f$ on $\partial \D$ is contained in $S^2$. If $S^2$ would be a totally real manifold, Gromov's compactness theorem would imply that $f$ is smooth on $\overline\D$. However, the presence of hyperbolic points requires additional considerations. The following assertion is due to \cite{hi}.

\begin{proposition}
The map $f$ extends continuously on $\overline\D$.
\end{proposition} 
\proof  Since the area of $f$ is bounded, for every point $p \in S^2$ and
$\varepsilon > 0$ small enough the intersection $f(\D) \cap B(p,\varepsilon)$
admits a finite number of connected disc-components in view of Lemma \ref{disc-connect}. We notice  here that if $p$ is a good hyperbolic point, then it follows from the description of the boundary behaviour of $f$ near  such a point established in Section 2 that the area of every disc in $f(\D) \cap B(p,\varepsilon)$ is separated from zero.

Let $p \in S^2_*$ be a totally real point of $S^2$. Lemma \ref{disc-connect}
may be applied to every map $f_k$ on $B(p,\varepsilon)$. Since the areas of
$f_k$ are bounded from above, the number of discs in $f_k(\D) \cap
B(p,\varepsilon)$ is bounded from above, independently of $k$. We use here the well-known fact (see for instance \cite{mc-sa}) that the areas of these discs are separated from zero. Applying to these discs Gromov's
compactness theorem, we obtain that their boundaries converge in
$B(p,\varepsilon) \cap \partial M$ to a continuous curve which is a connected
finite union  of smooth curves (as we will prove below, the bubbles do not
arise here so in fact the limit will be a single smooth curve). Thus, the
boundary of $f(\D)$ near a totally real point $p$ is a continuous curve. If
$p$ is a hyperbolic point, the boundary behavior of $f$ is described in
Section 2;  as it is proved there (see the basic representation (\ref{hyperb2})),  the boundary of $f(\D)$ near such a point also is a continuous curve. 

Now we may proceed the proof of the continuity of $f$. This is a slight modification of the classical Geometric Function Theory argument.

Suppose by contradiction that $f$ does not extend continuously to a point $\zeta_0 \in \partial \D$. Then there exist 
two sequences $(\zeta_n)$, $(\tilde \zeta_n)$ converging to $\zeta_0$ such
that $p_n = f(\zeta_n)$ and $\tilde p_n = f(\tilde \zeta_n)$ converge
respectively to $p_\infty \neq \tilde p_\infty$. Since the almost complex
structure $J$ is tamed by the symplectic form $\omega$,  they define
canonically a Riemannian metric $g$ (see, for instance, \cite{mc-sa}):

$$g(u,v) = (1/2)(\omega(u,Jv) + \omega(v,Ju))$$

 We will measure the distances and norms with respect to this metric (in fact, any Riemannian metric satisfying $g(X,X) \leq \omega(X,JX)$ is adapted). Let $d = dist(p_\infty,\tilde p_\infty)$. We may assume that $d(p_n,p_\infty) \leq d/3$, $d(\tilde p_n, \tilde p_\infty) \leq d/3$ for all $n$. Denote by $B(p,\varepsilon)$ the ball of radius $\varepsilon$ centered at $p$.

In view of the above description of 
the boundary behaviour of $f(\D)$ near $C(f,\partial \D)$, there exist maps (paths) $\lambda$ and $\tilde \lambda$, continuous on $[0,1]$, such that
$\lambda([0,1)) \subset f(\D) \cap B(p_\infty,d/3)$, $\tilde\lambda([0,1)) \subset f(\D) \cap B(\tilde p_\infty,d/3)$, and $\lambda(1) = p_\infty$, $\tilde\lambda(1) = \tilde p_\infty$. Furthermore, 
there exist increasing sequences $t_n,\tilde t_n \to 1$ such that $\lambda(t_n) = p_n$, $\tilde\lambda(t_n) = \tilde p_n$ (passing to a subsequence if necessarily).

Set $\Lambda(t) = f^{-1}(\lambda(t))$, $\tilde\Lambda(t) = f^{-1}(\tilde\lambda(t))$. Let $t_1(r)$ 
denote the smallest $t \in [0,1]$ such that $\Lambda(t) \in \partial B(\zeta_0,r)$. The
function $r \mapsto t_1(r)$ is decreasing, so is continuous except at most on a countable set of points of discontinuity. Put $\Lambda(t_1(r)) = \zeta(r) =
\zeta_0 + re^{i\tau(r)}$. Then the functions $r \mapsto \zeta(r)$ and $r
\mapsto \tau(r)$ also are continuous except on a countable set. We may define similarly the functions $\tilde \zeta$ and $\tilde \tau$.

Fix $r_0 > 0$ small enough. Using the polar coordinates $\zeta = re^{i\theta}$
on the disc $\zeta_0 + r_0\D$ and  integrating along the arc $[\tau(r), \tilde \tau(r)]$ of $\zeta_0 + r_0\partial\D$, we obtain :

\begin{eqnarray*}
d/3 \leq d(f(p_n),f(\tilde p_n)) \leq \int_{[\tau(r),\tilde\tau(r)]}
\left\vert \left \vert Df\left ( \frac{1}{r}\frac{\partial}{\partial
        \theta}\right ) \right \vert \right \vert r d\theta.
\end{eqnarray*}

By the Cauchy-Schwarz inequality we get :
\begin{eqnarray*}
(d/3)^2 & \leq & r^2\vert\tau(r) - \tilde \tau(r)\vert
\int_{[\tau(r),\tilde\tau(r)]} \left  \vert \left\vert Df\left (
      \frac{1}{r}\frac{\partial}{\partial \theta}\right ) \right \vert \right \vert^2  d\theta \\
&  \leq & 2\pi r^2 \int_{[\tau(r),\tilde\tau(r)]} \left\vert \left \vert
    Df\left ( \frac{1}{r}\frac{\partial}{\partial \theta}\right ) \right \vert
\right \vert ^2  d\theta.
\end{eqnarray*}

Let $dm(\zeta)$ denote the standard Lebesgue measure on $\C$ with respect to the variable $\zeta$. Dividing by $r$ and integrating with respect to $r$ from $\varepsilon$ to $r_0$ we obtain :

\begin{eqnarray*}
(d/3)^2 \ln(r_0/\varepsilon) & \leq & 2\pi \int_{[\varepsilon,r_0]} r dr
\int_{[\tau(r),\tilde\tau(r)]} \left\vert \left \vert Df\left (
      \frac{1}{r}\frac{\partial}{\partial \theta}\right ) \right \vert \right \vert^2  d\theta \\
& \leq & 2 \pi \int_{\D \cap (\zeta_0 + r_0 \D)}  \left\vert \left \vert
    Df\left ( \frac{1}{r}\frac{\partial}{\partial \theta}\right ) \right \vert
  \right \vert^2  dm(\zeta) \\
&  \leq & 2\pi \int_\D \omega \left ( Df\left ( \frac{1}{r}\frac{\partial}{\partial \theta}\right ),Df\left (i \frac{1}{r}\frac{\partial}{\partial \theta}\right ) \right) dm(\zeta) \\
& = & 2\pi area(f).
\end{eqnarray*}
Here in the last inequality we used the classical identity connecting the
symplectic area of a $J$-holomorphic curve and the $g$-norm of its
tangent vectors, see \cite{mc-sa}, pp. 20-21. Since $\varepsilon \to 0$, we obtain a contradiction which proves the Proposition. Q.E.D.

\subsection{Disc-bubbles do not arise at totally real points.}   The next step is the following

\begin{lemma} If $p \in f(\partial \D)$ is a totally real point of $S^2$, there are no disc bubbles arising at $p$.
\end{lemma}
\proof We could proceed similarly to \cite{ga-su}, but in the presence of
hyperbolic points it requires some global analysis of the characteristic foliation. So we use another argument due to \cite{ye}, easy to localize. Suppose that a disc-bubble $g$  arises at a totally real point $p = f_\infty(\zeta_1)$, $\zeta_1 \in \partial\D$. Let $\zeta_2 \in \partial \D$ be a point with $g(\zeta_2) = p$. By Proposition 2.6 of \cite{ga-su} the maps $f_\infty\vert_{ \partial \D}$ and $g \vert_{\partial \D}$ are embeddings near $\zeta_1$ and $\zeta_2$  respectively. Then we may fix open arcs $\gamma_j \subset \partial \D$ containing $\zeta_j$, $j=1,2$, such that $g(\gamma_2)$ lies on one side of $f_\infty(\gamma_1)$; in particular, $f_\infty(\gamma_1)$ and $g(\gamma_2)$ are tangent at $p$ (if not, the curves $f_k(\partial \D)$ would have self-intersections for large $k$, but they are embedded). Furthermore, because of the above mentioned relation for homological classes, the circles 
$f_\infty(\partial\D)$ and $g(\partial\D)$ have opposite orientations. Let
$\tau$ be a continuous unit tangent vector on the unit circle
$\partial\D$. Set $Y_k = \frac{\partial f_k}{\partial \tau} \left \vert \left
    \vert\frac{\partial f_k}{\partial \tau} \right \vert\right
\vert^{-1}$. Then, see \cite{ye}, it follows from the above description of a
disc-bubble  that there exists sequences $(\xi_k)$ and $(\tilde \xi_k)$ on
$\gamma_1$, converging to $\zeta_1$  and such that $f_k(\xi_k) \to p$,
$f(\tilde \xi_k) \to p$ and 
$$\lim_{k \to \infty} Y_k(\xi_k) = v, \lim_{k \to \infty} Y_k(\tilde \xi_k) = - v$$
where $v$ is a unit tangent vector to $f_\infty(\partial \D)$ at $p$. Let now $\chi$ be the characteristic foliation on $S^2$ , see \cite{ga-su}. Choose a continuous unit tangent vector field $X$ on 
 $S^2_*$  which is everywhere orthogonal to $\chi$ (with respect to the inner product $\bullet$ induced by some Riemannian metric). By Proposition 2.6 of \cite{ga-su} $f_\infty$ is transverse to $\chi$ at $p$; hence $v \bullet X(p) \neq 0$. Then for large $k$ the products $\frac{\partial f_k}{\partial \tau}(\xi_k)  \bullet X(f_k(\xi_k))$ and 
$\frac{\partial f_k}{\partial \tau}(\tilde \xi_k)  \bullet X(f_k(\tilde \xi_k))$ have opposite signs. By the intermediate value theorem, there exists a point $\eta_k \in \gamma_1$ such that 
$\frac{\partial f_k}{\partial \tau}(\eta_k)  \bullet X(f_k(\eta_k)) = 0$. This means that $f_k$ is tangent to $\chi$ at the point $f(\eta_k)$, which contradicts Proposition 2.6 of \cite{ga-su}. Q.E.D.

In particular, it follows by Gromov's compactness theorem that the sequence
$(f_k)$ converges to $f$ up to the boundary in each $C^k$-norm near every
totally real point $p \in S^2_*$. Furthermore, since $\partial M$ is $J$-convex, it follows by \cite{ga-su} that $f$ is transverse to $\partial M$ at $p$.

Now we study the behaviour of $(f_k)$ near a hyperbolic point.

\subsection{Dynamics of $(f^t)$ near a good hyperbolic point.}

Let $f^t$, $t \in \R$, be a 1-parameter family of embedded Bishop discs attached to
the totally real part $S^2_*$ of the sphere $S^2$ and converging to a non-constant
hyperbolic disc $f^\infty$. We call such a family {\it maximal}. Let $p$ be  a
good hyperbolic point in the boundary of $f_\infty$. We suppose that local
coordinates near $p$ are given by Definition \ref{def-good}. Therefore it follows from Section 2 that every disc $f^t$ is the graph of a function $g^t$ holomorphic in a domain $D^t$ in $\C$, in a neighborhood of the origin. The boundaries of the discs $f^t$ are disjoint so we may assume that the family $D^t$ of domains in $\C$ is either increasing or decreasing.

\begin{definition}
The family $f^t$ {\it  approaches $f^\infty$ from inside at $p$} (resp. from {\it outside}) if the family $D^t$ is increasing (resp. decreasing).
\end{definition}

In order to determine which case actually occurs near a given hyperbolic point, it is convenient to use the orientability of $S^2$. As above, consider  a 1-parameter family  $(f^t)$, $t \in \R$,  of Bishop discs attached to
the totally real part $S^2_*$. We observe that such a family $(f^t)$ provides $S^2$ with an orientation. This orientation is defined by pushing forward the form $dt \wedge d\theta$ via the map $$(t,\theta) \mapsto f^t(e^{i\theta}).$$ 

\begin{definition}
The family $f^t$ is called {\it positive} if this orientation coincides with an orientation already fixed on $S^2$. Otherwise a  maximal family is called {\it negative}.
 \end{definition}

We have the following
\begin{lemma}\label{posi-lemma}
A positive elliptic point generates a positive family of Bishop discs.
\end{lemma}
\proof Consider first the model case where $S^2$ is defined near an elliptic point
by the equation (\ref{hyperb1}) 
with $0 < \gamma < 1$. The Bishop discs near the origin are described in
\cite{ga-su} and are obtained by the intersection with the hyperplanes $\{ \Re w = t \}$. One readily sees that in this case the statement of lemma~\ref{posi-lemma} holds. Since the general case is a small perturbation of this model situation (see \cite{ga-su}) the assertion remains true. Q.E.D.

\begin{lemma}
\label{HypDiscLemma1}
Suppose that the family of discs $(f^t(\overline \D))_k$ approaches a good hyperbolic point $p$ from inside. Then (after a suitable reparametrization by Mobius transformations) the family $(f^t)$ converges (passing to a subsequence)  in $C(\overline\D)$ to a hyperbolic disc.
\end{lemma}
\proof We use local coordinates $(z,w)$ near $p$ as in Section 2, identifying $p$ with the origin. Then as $t \to \infty$ every disc $f^t(\D)$ near $0$
is the graph $\{ w = g^t(z) \}$ of a function $g^t$ holomorphic over a domain
$D^t$ with smooth boundary. Let $t_k \to \infty $ be an increasing  sequence
of values of parameter, $f_k = f^{t_k}$. The family $(D^t)$ is monotone  and
converges to a domain $D$ with $\partial D$ smooth off the
origin $D$; this  domain $D$ satisfies  Definition \ref{GoodApproach} of
Section 2.  Since the family $(D^t)$ is increasing, the limit domain 
$D$ can fill only one admissible region i.e. is asymptotic either to the
domain $D^+$ or $D^-$ from Definition \ref{GoodApproach}. 

 Next, the family $(g^{t_k})$ converges to a function $g$
holomorphic on $D$ near the origin. According to Section 2, $g$ is of class
$C^2(\overline D)$ and also satisfies Definition \ref{GoodApproach}.  As we saw
previously, the intersection of the limit nodal curve with a neighborhood of
the origin is a finite number of holomorphic discs in $M$ which are graphs
over $D$. But in our case the limit of discs
$(f_k(\D))$ is the graph of $g$ over $D$
near the origin so  only one single disc, the graph of $g$, appears in the
limit. Hence the origin is not a nodal point and  no disc bubble arises at the origin. If $t_n \to \infty$
is another sequence of parameters and $f_n = f^{t_n}$, we repeat the same
argument. Since $(D^t)$ is a monotone family, the domains $(D^{t_n})$ converge
to the same domain $D$ which is the limit of $(D^{t_k})$. Suppose that the
limit of $(f_n(\D))$ is the graph of $\tilde g$ over $D$ in a  neighborhood
$U$ of the origin. Then $\tilde g(z) = g(z) = \psi(z)$ for $z \in U
\cap \partial D$, where $\psi$ comes from Equation (\ref{hyperb1}). By the boundary
uniqueness theorem for holomorphic functions, $\tilde g = g$ on $D$. This means
that the family $(f^t)$ converges to a single disc which is the graph of $g$
and the absence of bubbles implies the convergence in $C(\overline\D)$. As it was mentioned  previously, near totally real points the convergence will be in every $C^k$ norm up to the boundary.  Q.E.D.

In the case of outside approach bubbles necessarily arise, but their structure is quite simple in view of a local description of hyperbolic discs from Section 2.

\begin{lemma}
\label{HypDiscLemma2}
Suppose that the family of discs $(f^t(\overline \D))_k$ approaches a good hyperbolic point $p$ from outside. Then near $p$ it converges in the Hausdorff distance to the union of two hyperbolic discs approaching $p$ from opposite sides in admissible regions described  in Section 2. Thus, the family of maps $(f^t)$ converges to one of these discs after a suitable reparametrization by Mobius transformations and the second disc may be viewed as a bubble.
\end{lemma}

\proof We use the characteristic foliation $\chi$ induced by $\partial M$ on $S^2$, see \cite{ga-su}. Since $p$ is a good hyperbolic point, we may assume that $S^2$ has the form (\ref{hyperb1}) near $p$ and $p = 0$ in these coordinates. Consider the projection $\pi:(z,w) \mapsto z$. Then the images  $\pi(\chi)$ foliate a neighborhood of the origin in $\C$. They are trajectories of a first order dynamical system 
$$
\left\{
\begin{array}{lll}
\Re \dot z  & = & (2\gamma+1) \Re z +  \alpha_1(2\gamma-1) \Im z  + O(z^2)\\
\Im \dot z & = & \alpha_1(2\gamma +1)  \Re z -  (2\gamma-1)\Im z + O(z^2)
\end{array}
\right.
$$
with $\alpha_1$, $\gamma$ given by the expansion of the local defining function
$\rho$ in Subsection 2.1. Then the origin is a saddle
point for this dynamical system. There are four trajectories through the
origin  tangent to   two lines through the origin; these lines are determined
by the eigenvalues of the linear part of the above system. They divide a
neighborhood of the origin in four regions, say, $\Omega_j$, $j=1,2,3,4$, which
are filled by other leaves of $\chi$ precisely as in the classical phase
portrait of a dynamical system near a saddle point. Now suppose by
contradiction that $(f^t)$ is a  family approaching $p$ from outside and
converging to a single disc. Then the limit disc is the graph over an
admissible approach region $D$, say asymptotic to the domain $ D^+$  from
Definition \ref{GoodApproach}. The discs $f^t(\D)$ are the graphs over domains
$D^t$ with smooth boundaries near the origin; these domains decrease to
$D$. But then for $t$ large enough their boundaries $\partial D^t$ intersect
at least three regions  $\Omega_j$ and one can find one region, say
$\Omega_1$, where the intersections $\partial D^t \cap \Omega_1$ form a
sequence of curves closed in $\Omega_1$, converging to the origin. But then,
since the origin is a saddle point, for a given $t$ one can find a leaf of the above dynamical system which is tangent to $\partial D^t$ at some point $a \in \Omega_1$. Therefore, there exists a leaf 
of the characteristic foliation $\chi$ tangent to $f^t(\overline\D)$ at some point. However, it is shown in \cite{ga-su} that this is impossible. Q.E.D.

Next, we have the following useful

\begin{lemma}
\label{positive-lemma}
If a positive (resp. negative) family $(f^t)$ approaches a positive hyperbolic point $p$, then it must approach from inside (resp. outside).
\end{lemma}
Once again, for the proof it suffices to consider the case where $S^2$ is
defined near a positive hyperbolic point by (\ref{hyperb1}) with $\gamma >
1$. In view of Lemma \ref{GoodApproach-lemma} the dynamics of the family $f^t$
near the origin is the same as the behaviour of the sections by the real
hyperplanes $\{ Re w = t \}$. We conclude the proof by checking the orientations of boundaries of discs corresponding to these sections. Q.E.D.

\vskip 0,2cm
Thus for every   maximal family $(f^t)$ two cases may occur. In the first case
(inside approach) 
the family $(f^t(\D))$ converges to the image of a single hyperbolic disc $f_1(\D)$ having a good
approach. Then necessarily we deal with the inside approach and the boundaries
of $f^t$ fill an approach region near a good hyperbolic point $p$. In the second case (outside
approach) the sequence $(f^t(\D))$ of discs 
converges to the union of images of two hyperbolic discs $f_1(\D)$ and $f_2(\D)$.

\section{Deformation of hyperbolic discs}  Recall again  \cite{ga-su}  that a Bishop disc with boundary glued to $S^2_*$ belongs to a 1-parameter family of Bishop discs with boundaries foliating an open piece of $S^2_*$. We establish here an analog for hyperbolic discs.  

\subsection{Gluing two hyperbolic discs into a single disc.} The main technical result here is the following

\begin{proposition}
\label{cut-off} Let $f_1$ and $f_2$ be two distinct hyperbolic discs at a good hyperbolic
    point $p$. Then given
    $\varepsilon > 0$ there exists an almost complex structure $J_\varepsilon$  and a sphere $S^2_\varepsilon$ with the following properties:
\begin{itemize}
\item[(a)] the structure $J_\varepsilon$ is   integrable near $p$ and
  coincides with $J$ outside a neighborhood of $p$; the sphere
  $S^2_\varepsilon$  also coincides with $S^2$ outside a neighborhood of $p$;
\item[(b)] $J_\varepsilon \to J$ in the $C^1$-norm  and $S^2_\varepsilon \to S^2$ in the $C^2$-norm as $\varepsilon \to 0$;
 \item[(c)] there exists a $J_\varepsilon$-holomorphic disc $f_\varepsilon$, $\varepsilon$-close to $f_1(\D)
    \cup f_2(\D)$ in the Hausdorff distance, coinciding with $f_1(\D)
    \cup f_2(\D)$ outside a neighborhood of $p$ and such that its boundary is glued to the totally real (with respect to $J_\varepsilon$) part  of
    $S^2_\varepsilon$. The family $(f_\varepsilon)$ tends to $f_1(\D) \cup f_2(\D)$ from outside as $\varepsilon \to 0$.
\end{itemize}
\end{proposition}
\proof {\it Step 1: smooth gluing of discs.} We assume that $S^2$ is given by (\ref{hyperb1}) near the origin which
is a good hyperbolic point.  Fix $\varepsilon > 0$ small enough.  According to Section 2, the discs $f_j$ have a good approach at $p$ and in particular, they approach $p$ from opposite regions. According to Lemma~\ref{GoodApproach-lemma}, the discs $f_j$ are the graphs of holomorphic (with respect to $J_{st}$) 
functions $w = g_j(z)$  over domains $D_j$ in $\C$ asymptotic to the origin
(the domains $D_j$ are asymptotic to   the domains $D^{\pm}$ by Definition \ref{GoodApproach}). Furthermore, the  expansions of $g_1$ and $g_2$ coincide at the origin up to the second order. In particular there exists a real number $\alpha$, $ 0 < \alpha < 1$, and a  $C^{2,\alpha}$-smooth real surface $\Pi$ in a neighborhood of the origin such that $f_j(\D)$ are contained in $\Pi$. Then $\Pi = \{ w = o(\vert z \vert^2) \}$. There exists 
 a  $C^{2,\alpha}$-coordinate diffeomorphism $\Phi_\varepsilon$ in a neighborhood $U$ of the origin with the following properties:
\begin{itemize}
\item[(i)]   $\Phi_\varepsilon$ coincides  with the identity map up to the second order at the origin, the restriction $\Phi_\varepsilon\vert_{f_j(\D)}$  is holomorphic with respect to $J_{st}$ and $\Phi_\varepsilon(f_j(\D)) = D_j$ , $j= 1,2$ ; 
\item[(ii)] in the new coordinates one has $\Pi = \{ w =  0 \}$ and $S^2$ has the form (\ref{hyperb0});
\item[(iii)] $\parallel \Phi_\varepsilon - id \parallel_{C^2} < \varepsilon$, where $id$ denotes the identity map.
\end{itemize}

{\it Step 2: local deformation of the structure.} Thus in the new coordinates one may identify the disc $f_j(\D)$ with the
domain $D_j$ in the axis $\{ (z,0) \}$ and there exists a small perturbation of
 $S^2$ satisfying (i), (ii)  and such that it has the form (\ref{hyperb1}); for
 simplicity of notations we still denote it by $S^2$. The structure $\tilde
 J_\varepsilon:= (\Phi_\varepsilon)_*(J)$ coincides with $J_{st}$ at the
 origin up to the first order.  Furthermore, $\tilde J_\varepsilon \vert_{D_j}
 = J_{st}$. Let $\tilde A_\varepsilon:= A_{\tilde J_\varepsilon}$ be the
 matrix of the deformation tensor of $\tilde J_\varepsilon$, see \cite{ga-su}
 (the complex matrix of the structure $\tilde J_\varepsilon$ in the
 terminology of \cite{su-tu3}). Recall that there exists a one-to-one correspondence between an almost complex structure and its deformation tensor, see \cite{su-tu3}.
In our case  $A_{\tilde J_\varepsilon}$ vanishes at the origin together with all first order partial derivatives and vanishes on the domains $D_j$.

{\it Step 3: local deformation of the sphere.} Let $V$ be a neighborhood of the origin in $\C$ and $\psi:V \to \R_+$ be  a
smooth function, $\psi(0) > 0$, with support compactly contained in $V$.
Let  $0 < \delta = \delta(\varepsilon) < \varepsilon$ be small enough.  Consider the surface $S^2_\varepsilon$ defined by 
\begin{eqnarray}
\label{SpPer1}
w =   z \overline z + \gamma \Re z^2  - \delta \psi(z).
\end{eqnarray}
Then $S^2_\varepsilon$ coincides with $S^2$ outside a small neighborhood of
the origin. The real surface  $\Pi = \{ w = 0 \}$ is $J_{st}$ -complex. Its intersection with $S^2_\varepsilon$ is a real curve in the $J_{st}$-totally real part of $S^2_\varepsilon$ and coincides  with the boundaries of $D_j$ outside a neighborhood of the origin.This curve bounds a $J_{st}$-holomorphic disc $\tilde f_\varepsilon$ on $\Pi$ coinciding with $D_1 \cup D_2$ outside a neighborhood of the origin. The family $\tilde f_\varepsilon$ tends to $D_1 \cup D_2$ from outside as $\varepsilon \to 0$.
Fix a smooth function $\chi(t)$ with $\chi(t) = 0$ for $t < 1$ and $\chi(t) =
1$ for $t > 2$. Consider the matrix $\hat A_\varepsilon(z,w) = \chi((\vert z
\vert + \vert w \vert)/\delta_1)\tilde A_\varepsilon(Z)$. Here  $0 < \delta_1
< \varepsilon$ is small enough. Then $\hat A_\varepsilon$ tends to $\tilde
A_\varepsilon$ in the $C^1$ norm as $\delta_1 \to 0$ (recall that $\tilde A_\varepsilon$ vanishes at the origin up to the first order). Furthermore, $\hat
A_\varepsilon$ vanishes near the origin, coincides with $\tilde A_\varepsilon$
outside a neighborhood of the origin and in the sectors $D_j$ where $\tilde
A_\varepsilon$ vanishes. Fixing $\delta < <\delta_1$  we obtain that the disc
$\tilde f_\varepsilon$ is holomorphic with respect to the almost complex
structure $\hat J_\varepsilon$ defined by the matrix $\hat
A_\varepsilon$. Hence the disc $f_\varepsilon = (\Phi_\varepsilon)^{-1}(\tilde
f_\varepsilon)$ and the structure $J_\varepsilon:=
(\Phi_\varepsilon)^{-1}(\hat J_\varepsilon)$ satisfy the assertion of the
lemma. Q.E.D.

\bigskip

{\bf Remark.} In the above proof we slightly perturbed the sphere $S^2$  and
the almost complex structure $J$  near a good hyperbolic point. Since the boundary $\partial M$
is strictly Levi convex near every hyperbolic point, it remains strictly Levi convex
after this $C^1$-perturbation of $J$.

\vskip 0,1cm
We assume that the point $p$ is a positive hyperbolic point and that the two hyperbolic discs $f_1$ and $f_2$ are limits of  two maximal positive families $(f_1^t)$ and $(f_2^t)$  of Bishop discs with Maslov indices equal to $0$. Applying Lemma \ref{cut-off}
we obtain a disc $f_\varepsilon$ approaching the union $f_1(\D) \cup f_2(\D)$ from outside as $\varepsilon$ tends to $0$.

\begin{lemma}\label{maslov-lemma}
The Maslov index of the disc $f_\varepsilon$ is equal to $0$.
\end{lemma}

\proof Fix $t$ large enough and consider the real oriented surface $E$ which is an open
piece of $S^2$   bounded by the curves  $f^t_1(\partial \D)$, $f^t_2(\partial \D)$ and $f_\varepsilon(\partial \D)$. We choose on these curves the orientation induced by the orientation of $E$ and denote the obtained loops by $\gamma_1$, $\gamma_2$ and $\gamma_3$ respectively. Thus, the orientations of $\gamma_j$, $j=1,2$ are opposite to the orientations induced by the discs $f^t_j$, $j=1,2$ and the orientation of $\gamma_3$ coincides with the orientation induced by $f_\varepsilon$. The boundary of $E$ is totally real with respect to the almost
complex structure $J_\varepsilon$ from 
Lemma  \ref{cut-off} and $E$ contains one positive hyperbolic point $p$. So
$I_+(E) = -1$ and $I_-(E) = 0$. The Maslov indices
of $f^t_j$, $j=1,2$ are equal to $0$. Hence $I_E(f^t_j\vert\partial\D)=1$, $j=1,2$ and   the sum of the indices $I_E(\gamma_j)$ is equal to $-2$. Then by Proposition \ref{ProInd} we obtain
$I_E(\gamma_3)  = +1$ and  the
Maslov index of $f_\varepsilon$ is equal to $0$. Q.E.D.

 In Proposition
\ref{cut-off} and Lemma \ref{maslov-lemma} we moved two hyperbolic discs $f_1$
and $f_2$ to a single disc which is "above" the hyperbolic point creating a
family of discs  approaching the hyperbolic point from outside; this family
bifurcates in the  initial hyperbolic discs $f_j$, $j=1,2$. This construction
is useful when $f_j$ are obtained as limits of two families of Bishop discs
approaching a hyperbolic point from inside. Of course, one can reverse this
dynamics of families of Bishop discs and adapt it for the case of outside
approach. Since we do not need it here, we drop the details.

\subsection{Hyperbolic chains} We considered above the case where hyperbolic discs touch exactly one hyperbolic point. It is easy to see that this construction admits a generalisation to the case where discs touch several hyperbolic points. More precisely, let $f_j$, $j=1,...,k$, be hyperbolic discs and let 
${\mathcal D} = \cup_{j=1}^k f_j(\D)$ be their union. We will use the notation
$G_j = f_j(\D)$; we stress out that the discs $G_j$ forming a chain are closed.

Suppose that $\overline {\mathcal D}$ is connected. We know that every disc touches a good hyperbolic point with a good asymptotic approach and by the uniqueness principle, if two discs touch the same hyperbolic point, they approach it from opposite regions.
Such family of discs (a nodal curve) is called in \cite{be-kl} {\it a hyperbolic chain}. A
chain ${\mathcal D}$ is called {\it saturated} if for every hyperbolic point
in $\overline{\mathcal D}$ both approach regions are filled by discs from
that chain. 

\begin{lemma}
\label{saturation}
If a chain ${\mathcal D}$ containing $k$ discs is not saturated, then
it contains at least $k$ hyperbolic points. 
\end{lemma}
\proof We proceed by induction in $k$. Consider the case where $k = 2$. If ${\mathcal D}$ contains only one hyperbolic point, by Section 2 the discs
$G_j$, $j= 1,2$, fill opposite admissible regions for that point and the chain
is saturated. Suppose now that the assertion of the Lemma holds for every chain containing at most $k-1$ discs. 

Consider a chain ${\mathcal D}$ containing $k$ discs and  denote by $m$ the number
of its hyperbolic points. Since the chain is not saturated, there exists a
hyperbolic point for which just one admissible region is filled, say, by the
disc  $G_k$. Removing $G_k$, we obtain a new chain ${\mathcal D}'$ containing $k-1$
discs and $m-1$ hyperbolic points. This new chain can not be saturated since
the disc $G_k$ is different from the discs forming  ${\mathcal D}'$ and approaches
at least one point of the chain ${\mathcal D}'$. Applying the
induction assumption, we conclude. Q.E.D.

Furthermore, in general a chain $\overline {\mathcal D}$ can be non-simply connected. For
example, one may think about a "triangle" formed by three hyperbolic
discs intersecting in three "vertices" that are the hyperbolic points. Then a
closed path formed by the diameters of these discs with the ends at the
hyperbolic points, is not homotopic to a point  in ${\mathcal D}$.

\begin{lemma}
\label{ChainConnect}
A  non-simply connected chain consisting of $k$ discs   contains at least $k$
hyperbolic points.
\end{lemma}
\proof Again we proceed by induction in $k$. Consider the case $k=2$. If
${\mathcal D}$ contains only one hyperbolic point, it is formed by two discs
glued together at this point and it is simply connected. Suppose that the
assertion of the Lemma holds for all chains containing at most $k - 1$ discs. 

Consider a chain ${\mathcal D}$ formed by $k$ discs and containing $m$
hyperbolic points. If $m < k$, then at least one disc, say $G_k$, contains only one
hyperbolic point $p$. Remove that disc from the chain obtaining the new chain
${\mathcal D}'$.  Since by the previous Lemma
the chain ${\mathcal D}$ is saturated, the hyperbolic point $p$
remains in the chain ${\mathcal D}'$ and  belongs to a single disc,
say $G_{k-1}$. Slightly deforming the disc $G_{k-1}$ and the almost complex
structure $J$ near $p$, we obtain a new chain (for the deformed structure) that does not contain $p$; hence it  contains $m-1$ hyperbolic points and it is simply connected by
assumption. Hence ${\mathcal D}'$ also is simply connected. But ${\mathcal
  D}$ is obtained by gluing the disc $G_k$ at a single point $p$, so it is
simply connected. Q.E.D.

 In the previous Subsection we proved  the existence of a deformation for
 saturated simply connected hyperbolic chains containing two discs and one
 hyperbolic point. 
This construction immediately generalises to the case of saturated simply connected chains containing  $k$ discs. We give details in the next section.

\section{Filling spheres}
Now we prove Theorem \ref{hyperb-theorem}. We proceed  by induction on the number $N$ of hyperbolic
points in $S^2$.

The case $N = 0$  is treated in \cite{ga-su}.

Consider the model  case $N=1$ : $S^2$ contains exactly one
hyperbolic point $H$. This point  may  be assumed positive changing the
orientation of $S^2$ if necessary. Then $S^2$ has three elliptic points $E_j$,
$j=1,2,3$, and necessarily two of them are positive because $I_+(S^2) =
I_-(S^2) = 1$. Denote by  $E_1$ and $E_2$ the positive elliptic points and by
$E_3$ the  negative one. Let $(f^t_j)$ be the maximal positive families of
Bishop discs generated by $E_j$, $j=1,2,3$  According to
Lemma~\ref{positive-lemma},  the families $(f_j^t)$, $j=1,2$ end up into two hyperbolic discs $f_1$ and
$f_2$. We point out that the families $(f_j^t)$, $j=1,2$ can not approach the negative elliptic point $E_3$: in that case they would touch 
the discs from the family $(f_3^t)$ and so would coincide with these discs by the uniqueness principle from \cite{ga-su}, which is impossible.
 By  Lemma~\ref{GoodApproach-lemma}, $f_1(\D)$ and $f_2(\D)$
approach $H$ from  two opposite regions. Applying Proposition 
\ref{cut-off} we obtain a family of discs $(f_\varepsilon)$ approaching $H$ from
outside; they are holomorphic with respect to a structure $J_\varepsilon$
which is a small perturbation of $J$ near $H$. Let now ${\tilde f}^t_3$ be the family
of $J_\varepsilon$-holomorphic Bishop discs generated by $E_3$. By the
uniqueness principle from  \cite{ga-su} the disc $f_\varepsilon$ is necessarily contained in the family $({\tilde f}^t_3)$. Passing to the limit as $\varepsilon \to 0$, we conclude that $S^2$ (after a small generic perturbation) is filled by boundaries of $J$-holomorphic discs.

\vskip 0,1cm
 Consider the case $N \geq 2$. We have $I_+(S^2) = I_-(S^2) = 1$.  Let
 $E_1$,...,$E_d$ be the positive elliptic points. Consider the positive families of Bishop discs $(f^t_j)$ generated by $E_j$. 

{\it Case 1.} The maximal families $(f^t_j)$ touch only positive hyperbolic points. Every family $(f^t_j)$ approaches a positive hyperbolic point from inside and fills one approach region. Denote by $f_j$, $j=1,...,d$,  the limit hyperbolic disc for every family. We can regroup these discs to a finite family of disjoint hyperbolic chains ${\mathcal D_l}$.
Since the number of   positive hyperbolic points is less than $d$, one of the
chains, say ${\mathcal D}_1$, contains more discs than hyperbolic points. Let $n$
and $m$ be respectively the numbers of the  discs and the hyperbolic points in
${\mathcal D}_1$. Since $m < n$,  this chain is saturated and simply
connected.  Applying to this chain the deformation construction from the
previous section (Proposition \ref{cut-off}), we swept an open subset of $S^2$ containing  the hyperbolic
points from ${\mathcal D}_1$ by the boundaries of discs holomorphic with
respect to an $\varepsilon$-perturbed almost complex structure. We extend
these families past the hyperbolic points; since the chain saturated and
simply connected, we obtain   $(n-m)$ positive families of Bishop discs of the
Maslov index $0$. Thus we removed the same number of discs and hyperbolic
points and may proceed by induction. Finally we obtain a filling of $S^2$ by
boundaries of discs holomorphic with respect to an almost complex structure
$J_\varepsilon$ obtained from $J$ by an $\varepsilon$-perturbation 
near every hyperbolic point. Passing to the limit as $\varepsilon \to 0$, we conclude.

{\it Case 2.} There are negative hyperbolic points which are limits of the
families $(f^t_j)$. Once again we regroup the limit hyperbolic discs $f_j$
into disjoint chains ${\mathcal D}_j$. Since a positive family of discs approaches a negative
hyperbolic point from outside, it fills both approach regions by Section 4 and
the limit chains are saturated at every negative hyperbolic point;
furthermore, every negative hyperbolic point attracts precisely one positive
family of discs. Next we may
apply a deformation construction similarly to   Proposition \ref{cut-off} at every negative hyperbolic
point and slightly deform every chain near such a point; we also suitably
deform the almost complex structure  so that the new chains remain
holomorphic. The new chains contain only  positive hyperbolic points  and by
the previous argument, one of these chains  is saturated and simply
connected. Since such a chain is
obtained by a small deformation of an old one, we conclude that one of the initial chains,
say ${\mathcal D}_1$, is saturated and simply connected. Then we extend this
chain through positive hyperbolic points by Proposition  \ref{cut-off}
replacing two positive families and at least one positive hyperbolic point by a
family of discs of Maslov index $0$. If the above  chains contain only negative
hyperbolic points, we consider families of Bishop's discs starting from
negative hyperbolic points. Such a family fills both admissible regions near a
positive hyperbolic point. Hence in every chain obtained from negative families
of discs  every positive hyperbolic point attracts precisely  one family of
discs. Since $I(S) = 2$, at least one of these chains contains a negative
hyperbolic point and we apply the above argument extending the chain through
that point by Proposition \ref{cut-off}. Thus, in any case we remove at least
one hyperbolic point and two families of discs replacing them by a single
family. By the induction assumption, we conclude the proof of Theorem
\ref{hyperb-theorem}. 

Thus, we obtain a real hypersurface $\Gamma$ with boundary $S^2$. By
construction, $\Gamma$ is a real smooth hypersurface foliated by a 1-parameter
family of 
$J$-holomorphic discs; it is obvious from the above construction that $\Gamma$
is diffeomorphic to the 3-ball. 

\section{Concluding remarks}

In this section we compare Theorem \ref{hyperb-theorem} with related
results.We do not discuss numerous applications of this Theorem to
symplectic and contact geometry; see, for instance  \cite{el}.

{\bf 0.} A true breakthrough in the study of filling of 2-spheres in presence
of hyperbolic points was done by E.Bedford-W.Klingenberg \cite{be-kl}. They
consider the case of spheres with elliptic and good hyperbolic points
generically embedded into a strictly pseudoconvex hypersurface
in $\C^2$. This work remains an important reference in the subject.

{\bf 1.} In the interesting and important work by R.Hind \cite{hi} the following situation is considered.
Let $(M,J,\omega)$ denote a symplectic manifold with a tame almost complex structure containing no holomorphic spheres of negative self-intersection.
Let $\Omega$ be a smoothly bounded domain with Levi convex boundary $\partial\Omega$. Suppose that $\partial \Omega$ is not the cartesian product of a holomorphic sphere with the circle $S^1$ and let $S^2$ be a real 2-sphere generically embedded into $\partial\Omega$. Suppose that $J$ is integrable in a neighborhood of every hyperbolic point of $S^2$. Then, if necessarily after a $C^2$ perturbation  in a neighborhood of the complex points, there exists a filling of $S^2$ by boundaries of holomorphic discs. The work contains a detailed description of the properties of this filling. Admitting that
the result holds in the case where $\partial \Omega$ is strictly Levi convex, the author uses Y.Eliashberg-W.Thurston's theorem \cite{el-th} on approximation of Levi convex boundaries by strictly Levi convex ones. However, in the almost complex setting the corresponding result on filling in the strictly Levi convex case was never proved (it was announced by Y.Eliashberg in \cite{el}).  

The present paper fixes that gap. In fact, we prove much more since Theorem
\ref{hyperb-theorem} is established in the Levi convex case under the
assumption that there are no non-constant holomorphic discs in the
boundary. This is the main case considered in the work \cite{hi} since the
case where the boundary contains holomorphic discs can be reduced to the
previous one by approximation of the boundary, see \cite{hi}. We also point
out that many other technical simplifications of R.Hind's work are given in
\cite{ga-su} using exhaustion
plurisubharmonic functions.

Note that in general, the condition of the Levi convexity cannot be dropped. It was first observed by Y.Eliashberg \cite{el}; later J.Fornaess-D.Ma \cite{fo-ma} constructed an explicit example.

{\bf 2.}  N.Kruzhilin proved in \cite{kr} the existence of a filling of a two-sphere
generically embedded into a strictly Levi convex hypersurface in $\C^2$ (with
the standard complex structure) without any assumption that the hyperbolic points
are good. His description of the boundary behaviour of hyperbolic discs is
substantially more complicated. It is natural to think that a combination of
his techniques with the methods of the present work will allow to obtain
Theorem \ref{hyperb-theorem} without assuming that  the hyperbolic points
are good. Another (more general) open question is if Theorem
\ref{hyperb-theorem} remains true when the almost complex structure $J$ is not
supposed to be integrable near hyperbolic points. However, we must point out
that for many applications the condition of the integrability of $J$ near
hyperbolic points and the assumption that these points are good are not
restrictive and naturally hold. Indeed, by a small perturbation of an almost
complex structure a given two-sphere can be lead into a position where
the assumptions of Theorem \ref{hyperb-theorem} hold; the integrability of
the structure near complex points also often can be assumed. Thus, Theorem
\ref{hyperb-theorem} is in general sufficient for filling a single sphere.
So it is quite possible that the necessary technical
difficulties to answer the above open questions will not correspond to
the impact. 

{\bf 3.} As mentioned in \cite{ga-su}, the condition that
$(M,J,\omega)$ contains no non-constant  holomorphic spheres 
may be weakened. We do not develop this subject here referring to
the concluding remarks in \cite{ga-su}. 

{\bf 4.} Many parts of the proof presented in our article still  work in the case where instead of the real sphere $S^2$ one considers a compact real surface (with or without boundary) contained in a pseudoconvex hypersurface and  admitting a finite number of complex points.  Here  suitable assumptions on  the numbers  of elliptic and hyperbolic points are necessary. We hope that it will be useful in further applications.

\end{document}